\documentclass[12pt,a4paper]{amsart}
\usepackage{amssymb, amsmath}

\usepackage{fullpage}
\usepackage{mathrsfs}
\usepackage{enumerate}
\usepackage{amsthm}
\usepackage{amsmath}
\usepackage{amsfonts}
\usepackage{amssymb}

\usepackage[nobysame,initials,logical-quotes]{amsrefs}
\usepackage[ocgcolorlinks,hyperfootnotes=false,colorlinks=true,citecolor=blue,linkcolor=blue,urlcolor=blue]{hyperref}
\usepackage{graphicx}
\usepackage{color}
\usepackage{graphics}
\usepackage{eepic}
\usepackage{nicefrac}

\newtheorem{thm}{Theorem}[section]

\newtheorem{lemma}[thm]{Lemma}

\newtheorem{proposition}[thm]{Proposition}

\theoremstyle{definition}

\newtheorem{definition}[thm]{Definition}
\newtheorem{remark}[thm]{Remark}
\newtheorem{example}[thm]{Example}

\newcommand{\ignore}[1]{}

\usepackage[colorinlistoftodos]{todonotes}

\newcommand{\bydef}{\stackrel{\rm def}{=}}

\usepackage{mathtools}

\numberwithin{equation}{section}

\def\textmatrix#1&#2\\#3&#4\\{\bigl({#1 \atop #3}\ {#2 \atop #4}\bigr)}
\def\dispmatrix#1&#2\\#3&#4\\{\left({#1 \atop #3}\ {#2 \atop #4}\right)}
\numberwithin{equation}{section}

\begin{document}
 \title{A note on the dilation of a certain family of tetrablock contractions\footnote{MSC 2020: 47A20, 47A13}}
\author{Tirthankar Bhattacharyya}
\address{Department of Mathematics\\
Indian Institute of Science\\
Bangalore 560012, India}
\email{tirtha@iisc.ac.in}
\author{Mainak Bhowmik}
\address{Department of Mathematics\\
Indian Institute of Science\\
Bangalore 560012, India}
\email{mainakb@iisc.ac.in}
\maketitle

\begin{abstract}
We find an explicit tetrablock isometric dilation for every member $(A_\alpha, B, P)$ of a family of tetrablock contractions indexed by a parameter $\alpha$ in the closed unit disc (only the first operator of the tetrablock contraction depends on the parameter). The dilation space is the same for any member of the family. Explicit dilation for the adjoint tetrablock contraction $(A_\alpha^*, B^*, P^*)$ for every member of the family mentioned above is constructed as well. This example is important because it has been claimed in the literature that this example does not have a dilation. Taking cue from this construction and using Toeplitz operators on $H^2_{\mathbb D}(\mathcal D_P)$, we obtain necessary and sufficient conditions for a tetrablock contraction to have a certain type of tetrablock isometric dilation.
\end{abstract}

\section{Introduction}\label{intro}
Given a compact subset $K$ of $\mathbb{C}^d$ for any $d=1,2,3,\dots$ and a tuple $\underline{T} = \left(T_1,T_2,\dots,T_d \right)$ of commuting bounded operators on a Hilbert space $\mathcal{H}$, the set $K$ is called a \textit{spectral set} for $\underline{T}$ if the Taylor joint spectrum, $\sigma(\underline{T})\subseteq K $ and $$ \lVert f(\underline{T}) \rVert \leq \lVert f \rVert $$ for any rational function $f:\mathbb{C}^d \to \mathbb{C}$ with poles off $K$, where $$ \lVert f \rVert = \sup\left\lbrace |f(\underline{z})| : \underline{z} \in K \right\rbrace  $$ and $f(\underline{T})$ is formed according to Taylor's functional calculus.\\
The set $K$ is called a \textit{complete spectral set} for $\underline{T}$ if given any matrix valued rational function, $f= \left(\left( f_{ij}\right)\right): \mathbb{C}^d \to \mathbb{M}_n(\mathbb{C})$ with poles off $K$, for $n=1,2,\dots $ we have $$ \lVert f(\underline{T})\rVert \leq \lVert f \rVert =\sup\left\lbrace \lVert \left((f_{ij}(\underline{z}))\right)\rVert_{op} : \underline{z} \in K \right\rbrace $$
where $$ f(\underline{T})= \left(\left( f_{ij}(\underline{T}) \right)\right)_{i,j=1}^n $$ is a bounded operator from the direct sum of $n$ copies of $\mathcal{H}$ to itself.
\\
\begin{definition}\label{dilation-def}
A $\partial K $-\textit{normal dilation} for a tuple $\underline{T}$ with $\sigma(\underline{T}) \subseteq K $ is a tuple of commuting bounded normal operators $\underline{N} = \left(N_1,N_2,\dots,N_d \right)$ on a Hilbert space $\mathcal{K}$ containing $\mathcal{H}$ as a closed subspace satisfying $$ f(\underline{T}) = \text{P}_{\mathcal{H}} f(\underline{N}) \vert_{\mathcal{H}}  $$
for every rational function $f$ on $\mathbb{C}^d$ with poles off $K$ and $\sigma(\underline{N})\subseteq \partial K$, the Silov boundary of $K$ with respect to the algebra of functions which are continuous on $K$ and holomorphic in the interior of $K$.
\end{definition}
\begin{thm}{\textbf{(Arveson)}}
Let $\underline{T}$ have $K$ as a spectral set. Then $\underline{T}$ has $K$ as a complete spectral set if and only if $\underline{T}$ has a $\partial K$-normal dilation.
\end{thm}

We say that rational dilation \textit{holds} on $K$ if $K$ is a complete spectral set for $\underline{T}$ whenever $K$ is a spectral set for $\underline{T}$. Else, rational dilation \textit{fails}. The following cases are well-known. The rational dilation
\begin{enumerate}
\item[1.] holds on the unit disc $\overline{\mathbb{D}} = \left\lbrace z\in \mathbb{C} : |z|\leq 1 \right\rbrace$ (von Neumann \cite{vN}),
\item[2.] holds on the bidisc $\overline{\mathbb{D}}^2$ (Ando \cite{Ando}),
\item[3.] fails in general on polydisc $\overline{\mathbb{D}}^d$ for $d\ge 3$ (Varopoulos \cite{Var}),
\item[4.] holds in an annulus (Agler \cite{agler}),
\item[5.] fails in general on a planar domain with two or more holes (Agler-Harland-Rafael \cite{ahr} and Dritschel-McCullough \cite{DM}),
\item[6.] holds on the symmetrized bidisc (Agler-Young \cite{ay-jfa}),
\item[7.] fails on any norm unit ball in $\mathbb{C}^d$ for $d\geq 3$ (Paulsen \cite{Paulse-JFA} and Pisier \cite{Pisier-os}).
\end{enumerate}

In the backdrop of this, we consider a polynomially convex and compact subset $\overline{\mathbb{E}}$ of $\mathbb{C}^3$ defined as follows. Let $\mathbb M_2$ denote the algebra of all $2 \times 2$ complex matrices equipped with the operator norm. Define $\pi : \mathbb M_2 \rightarrow \mathbb C^3$ by
$$ \pi \left(
         \begin{array}{cc}
           a_{11} &  a_{12} \\
           a_{21} & a_{22} \\
         \end{array}
       \right) = (a_{11}, a_{22}, \det A).$$

The closed $tetrablock$ $\overline{\mathbb{E}}$ is defined to be the set of all $\pi(A)$ such that $A$ is in $\mathbb{M}_2$ with $\lVert A \rVert \leq 1$. There has been an interesting history of the theory of commuting triples $(A,B,P)$ of operators with $\overline{\mathbb{E}}$ as a spectral set. Such a triple is called a \textit{tetrablock contraction} or an {\em $\overline{\mathbb{E}}$-contraction}. A special class of tetrablock contractions was dilated in \cite{Bhattacharyya}.

{\em Whether rational dilation holds on the tetrablock is an open question.}

In the context of the tetrablock, a boundary normal dilation consists of a triple $(N_1, N_2, N_3)$ of commuting normal operators with joint spectrum contained in the distinguished boundary
$$ b\mathbb E =\left\lbrace \pi(A) : A \in \mathbb{M}_2 \text{ is a unitary matrix}\right\rbrace.$$
Such a normal triple is called a {\em tetrablock unitary}. A  {\em tetrablock isometry} is the restriction of a tetrablock unitary to a joint invariant subspace. Although the Definition \ref{dilation-def} of dilation demands that given a tetrablock contraction $(A, B, P)$, one must construct a normal triple as above, it is straightforward (exactly as in the case of the unit disc, the bidisc or the symmetrized bidisc) to see that it is enough to construct a Hilbert space $\mathcal K$ containing $\mathcal H$ and a tetrablock isometry $(V_1, V_2, V_3)$ on $\mathcal K$ such that
$$V_1^*|_{\mathcal H} = A^*, V_2^*|_{\mathcal H} = B^* \text{ and } V_3^*|_{\mathcal H} = P^*.$$
The space $\mathcal H$ is then $co-invariant$ under the tetrablock isometry. This is what we do in this note for a family of $\overline{\mathbb{E}}$-contractions. The characterization of tetrablock isometries that we shall use is as follows:
\begin{thm}{(Theorem 5.7 in \cite{Bhattacharyya})} \label{thm-e-contraction}
A commuting triple $(V_1,V_2,V_3)$ of bounded operators on a Hilbert space is a tetrablock isometry if and only if
\begin{enumerate}
\item $V_3$ is an isometry,
\item $r(V_1)\leq 1$, $r(V_2)\leq 1$ (equivalently, $\| V_1\| \le 1$ and $\|V_2\| \le 1$) and
\item $V_1=V_2^*V_3$
\end{enumerate}
where $r$ stands for the spectral radius.

\end{thm}

Section \ref{sec-1} does the construction of dilation for a family motivated by the example from \cite{Pal}. Section \ref{adj} uses a general theorem about adjoints and ties with Section \ref{sec-1}.  Dilations for the adjoints of the tetrablock contractions considered in Section \ref{sec-1} are constructed in this section. Section \ref{DilThmSec} takes a cue from Section \ref{sec-1} and investigates those tetrablock contractions which have a special kind of dilation.

\section{Dilation of a family of tetrablock contractions}
\label{sec-1}
In this section, we construct an explicit tetrablock isometric dilation for the example found in Section 5 of \cite{Pal}. The same is done for the adjoint of this example in the next section. See also Ball and Sau \cite{Ball-Sau} who proved that the "necessary condition" in \cite{Pal} is not necessary, but did not conclude anything about existence or non-existence of a tetrablock isometric dilation for this example.

For a Hilbert space $\mathcal{E}$, the infinite direct sum $\mathcal{E} \oplus \mathcal{E} \oplus \cdots$ will be denoted by $\ell ^2(\mathcal{E})$. It will often be identified with $H^2_{\mathbb D}(\mathcal{E})$, the Hardy space of $\mathcal{E}$-valued functions on the open unit disc. Accordingly $T_z$ will denote the unilateral shift of multiplicity equal to the dimension of $\mathcal{E}$ on $H^2_{\mathbb D}(\mathcal{E})$ $\left(\text{or,}\, \ell^2(\mathcal{E}) \right)$, i.e.,
$$ T_zf(z) = zf(z) \text{ for } f \text{ in } H^2_{\mathbb D}(\mathcal{E}) \text{ or equivalently } T_z(e_0, e_1, \ldots ) = (0, e_0, \ldots ) \text{ for } e_n \text{ in } \mathcal E .$$
\begin{example}\label{Pal-example}
Consider $\mathcal{H}= \ell ^2(\mathbb{C}^2)\oplus \ell ^2(\mathbb{C}^2) \oplus \ell ^2(\mathbb{C}^2)\oplus \ell ^2(\mathbb{C}^2)$ . Let $A_\alpha,B,P$ be defined by the following block operator matrices on $\mathcal{H}$ with respect to the decomposition above.
$$
A_\alpha=\begin{bmatrix}
 0 & 0 & 0 & 0\\
 0 & 0 & 0 & 0\\
 0 & 0 & H & 0\\
 0 & 0 & 0 & 0
 \end{bmatrix}, \, B= \begin{bmatrix}
 0 & 0 & 0 & 0\\
 0 & 0 & 0 & 0\\
 0 & 0 & 0 & 0\\
 0 & 0 & 0 & 0
 \end{bmatrix},\, P = \begin{bmatrix}
 0 & 0 & 0 & 0\\
 0 & 0 & 0 & 0\\
 0 & T_z & 0 & 0 \\
 I & 0 & 0 & 0
 \end{bmatrix}
$$
where $H$ on $\ell^2(\mathbb{C}^2)$ is defined as follows:
$$
H(c_0,c_1,c_2,\cdots) := (H_1c_0, 0,0,\cdots) \,\, \text{for}\,\, (c_0,c_1,c_2,\cdots) \in \ell^2(\mathbb{C}^2)
$$
and $H_1$ is the $2 \times 2$ matrix $H_1 = \textmatrix 0 & \alpha \\ 0 & 0 \\ $ for $\alpha \in \overline{\mathbb{D}}$. We shall often suppress $\alpha$ and write the tetrablock contraction as $(A, B, P)$. It is understood that $\alpha$ is present in $A$ in the way above.

Note that the product of any two of $A,B,P$ is $0$.
We shall not prove $(A,B,P)$ to be a tetrablock contraction here because we shall dilate $(A,B,P)$ to a tetrablock isometry which will automatically prove it to be a tetrablock contraction.
 \end{example}
 {\em The dilation will be for any scalar  $\alpha \in \overline{\mathbb{D}}$. The example in section 5 of \cite{Pal} is for $\alpha =\frac{1}{4} $}. The defect operator of $P$ is
  $$D_P^2 = I - P^*P =\begin{bmatrix}
  I & 0 & 0 & 0\\
  0 & I & 0 & 0 \\
  0 & 0 & I & 0\\
  0 & 0 & 0 & I
 \end{bmatrix} - \begin{bmatrix}
0 & 0 & 0 & I\\
0 & 0 & T_z^* & 0\\
0 & 0 & 0 & 0\\
0 & 0 & 0 & 0
 \end{bmatrix}
 \begin{bmatrix}
 0 & 0 & 0 &0\\
 0 & 0 & 0 &0\\
 0 & T_z & 0 & 0\\
 I & 0 & 0 & 0
 \end{bmatrix}= \begin{bmatrix}
 0 & 0 & 0 &0\\
 0 & 0 & 0 &0\\
 0 & 0 & I & 0\\
 0 & 0 & 0 & I
 \end{bmatrix}
 $$
 So $ D_P$ is the projection onto
$$\mathcal{D}_P = Ran(D_P) = \{0\} \oplus \{0\} \oplus \ell ^2(\mathbb{C}^2)\oplus \ell ^2(\mathbb{C}^2).$$

Let \begin{align} \label{F1F2} F_1  = 0 \oplus \begin{bmatrix}
H & 0\\
0 & 0
\end{bmatrix}
\text{and} \, \, F_2 = 0  \end{align}
on $\mathcal{D}_P =\left( \{0\} \oplus \{0\} \right) \oplus \left( \ell ^2(\mathbb{C}^2)\oplus \ell ^2(\mathbb{C}^2)\right)$. These $F_1$ and $F_2$ satisfy $$ A-B^*P = A = D_P F_1 D_P , \,\,  B-A^*P = D_P F_2 D_P \,\, \text{and} \,\, F_1F_2 = F_2F_1 .$$ These can be seen by easy computations. We have $\left[ F_2 , F_2^* \right] = 0$ since $F_2 = 0$. However,  $$[F_1, F_1^*] = 0 \oplus \begin{bmatrix}
HH^*- H^*H & 0 \\
0 & 0
\end{bmatrix} \neq 0 $$
since $H$ is non-normal.

Since $H_1^2 = 0$ we have $H^2 = 0$. Therefore $F_1^2 = 0$. It is well-known that if a commuting triple $(A,B,P)$ has a tetrablock isometric dilation $(V_1,V_2,V_3)$ (which implies that $(A,B,P)$ is a tetrablock contraction) and $W_3$ is the minimal isometric dilation of $P$ then the solutions $F_1$ and $F_2$ of $A-B^*P= D_P F_1 D_P$ and $B-A^*P = D_P F_2 D_P$ have to satisfy, among other things, $\left[ F_1, F_1^*\right] = \left [ F_2, F_2^* \right]$ which is not the case here.
\subsection*{Construction of the dilation triple :}
Consider a triple of bounded operators $(V_1,V_2,V_3)$ on $\mathcal{K}= \mathcal{H}\oplus \mathcal{D}_P \oplus \mathcal{D}_P \oplus \cdots$ as follows:
$$
V_1= \left[ \begin{array}{c | c c c c c}
A & 0 & 0 & 0 & 0 &\cdots \\
\hline
F_1^*D_P & F_1 & 0 & 0 & 0 & \cdots \\
0 & F_1^* & F_1 & 0 & 0 &\cdots \\
0 & 0 & F_1^* & F_1 & 0 & \cdots \\
0 & 0 & 0 & F_1^* & F_1 & \cdots \\
\cdots & \cdots & \cdots & \cdots & \cdots & \cdots
\end{array} \right] ,\,\,
V_2 = \left[ \begin{array}{c | c c c c c}
B & 0 & 0 & 0 & 0 &\cdots \\
\hline
F_1D_P & 0 & 0 & 0 & 0 & \cdots \\
F_1^*D_P & F_1 & 0 & 0 & 0 & \cdots \\
0 & F_1^* & F_1 & 0 & 0 & \cdots \\
0 & 0 & F_1^* & F_1 & 0 & \cdots \\
\cdots & \cdots & \cdots & \cdots & \cdots & \cdots
\end{array}\right]
$$
$$ V_3 = \left[\begin{array}{c | c c c c c}
P & 0 & 0 & 0 & 0 & \cdots \\
\hline
0 & 0 & 0 & 0 & 0 & \cdots \\
D_P & 0 & 0 & 0 & 0 & \cdots \\
0 & I & 0 & 0 & 0 & \cdots \\
0 & 0 & I & 0 & 0 & \cdots \\
\cdots & \cdots & \cdots & \cdots & \cdots & \cdots
\end{array}\right].
$$
 We shall show that $(V_1,V_2,V_3)$ is a commuting triple. To that end, we compute
$$V_1V_2 = \left[\begin{matrix}
AB & 0 & 0 & \cdots \\

F_1^*D_PB + F_1^2D_P & 0 & 0 & \cdots \\

 F_1^* F_1D_P + F_1F_1^*D_P  &  F_1^2 & 0 & \cdots \\

 {F_1^*}^2D_P &  F_1^*F_1 + F_1F_1^* & F_1^2 & \cdots \\

0 & {F_1^*} ^2 & F_1^* F_1 + F_1F_1^* & \cdots \\

0 & 0 & {F_1^*}^2 & \cdots \\

\cdots & \cdots & \cdots & \cdots
\end{matrix} \right] $$
and
$$ V_2V_1 = \left[\begin{matrix}
BA & 0 & 0 &\cdots \\
F_1D_PA & 0 & 0 & \cdots \\
F_1^*D_PA + F_1F_1^*D_P & F_1^2 & 0 & \cdots \\
{F_1^*}^2 D_P & F_1^*F_1 + F_1F_1^* & F_1^2 & \cdots \\
0 & {F_1^*}^2 & F_1^*F_1 + F_1F_1^* & \cdots \\
\cdots & \cdots & \cdots & \cdots
\end{matrix}\right].
$$
Since $F_1^2 = 0$ and $B=0$, we have
$$
F_1^*D_PB + F_1^2D_P = 0 = \left(0 \oplus \begin{bmatrix}
H & 0\\
0 & 0
\end{bmatrix}\right) \left( 0 \oplus \begin{bmatrix}
I & 0\\
0 & I
\end{bmatrix}\right) \left(0 \oplus \begin{bmatrix}
H & 0\\
0 & 0
\end{bmatrix}\right) = F_1D_P A
$$
which proves that the (2,1)-th entries in $V_1V_2$ and $V_2V_1$ are equal.
Again $$
F_1^* F_1D_P + F_1F_1^*D_P = 0 \oplus \begin{bmatrix}
H^*H + HH^* & 0\\
0 & 0
\end{bmatrix} = F_1^*D_PA + F_1F_1^*D_P
$$
which shows that the (3,1)-th entries in  $V_1V_2$ and $V_2V_1$ are equal. Hence $V_1V_2 = V_2V_1$.\\
To check that $V_1V_3 = V_3V_1$, we compute
$$ V_1V_3 = \begin{bmatrix}
 AP & 0 & 0 & 0 & \cdots \\
 F_1^*D_P P & 0 & 0 & 0 & \cdots \\
 F_1D_P & 0 & 0 & 0 & \cdots \\
 F_1^*D_P & F_1 & 0 & 0 & \cdots \\
 0 & F_1^* & F_1 & 0 & \cdots \\
 0 & 0 & F_1^* & F_1 & \cdots \\
 \cdots & \cdots & \cdots & \cdots & \cdots
 \end{bmatrix} \,\,\text{and}\,\,
V_3V_1 =
 \begin{bmatrix}
PA & 0 & 0 & 0 & \cdots \\
0 & 0 & 0 & 0 & \cdots \\
D_PA & 0 & 0 & 0 & \cdots \\
F_1^*D_P & F_1 & 0 & 0 & \cdots \\
0 & F_1^* & F_1 & 0 & \cdots \\
0 & 0 & F_1^* & F_1 & \cdots \\
\cdots & \cdots & \cdots & \cdots & \cdots
\end{bmatrix}
.$$
But $ F_1^*D_P P = 0$ since $H^* T_z = 0 $. An easy computation shows $ D_P A  = F_1D_P$. Therefore, each entry of the block operator matrix $V_1 V_3$ is the same as the corresponding entry in the matrix of $V_3 V_1$. So, $V_1$ and $V_3$ commute. Similarly,$$ V_2V_3 =
\begin{bmatrix}
BP & 0 & 0 & \cdots \\
F_1D_PP & 0 & 0 & \cdots \\
F_1^*D_PP & 0 & 0 & \cdots \\
F_1D_P & 0 & 0 & \cdots \\
F_1^*D_P & F_1 & 0 & \cdots \\
0 & F_1^* & F_1 & \cdots \\
\cdots & \cdots & \cdots & \cdots
\end{bmatrix}  \,\,\text{and}\,\,
 V_3V_2 =
 \begin{bmatrix}
PB & 0 & 0 & \cdots \\
0 & 0 & 0 & \cdots \\
D_PB & 0 & 0 & \cdots \\
F_1D_P & 0 & 0 & \cdots \\
F_1^*D_P & F_1 & 0 & \cdots \\
0 & F_1^* & F_1 & \cdots \\
\cdots & \cdots & \cdots & \cdots
\end{bmatrix}.
$$We have $$ F_1D_P P = \begin{bmatrix}
0 & & &  \\
 & 0 & &  \\
 & & H & \\
  & & & 0
\end{bmatrix}
\begin{bmatrix}
0 & & & \\
 & 0 & & \\
  & & I & \\
  & & & I
\end{bmatrix}
\begin{bmatrix}
0 & 0 & 0 & 0 \\
0 & 0 & 0 & 0 \\
0 & T_z & 0 & 0 \\
I & 0 & 0 & 0
\end{bmatrix} = \begin{bmatrix}
0 & & & \\
 & 0 & & \\
  & &0 & \\
  & & & 0
\end{bmatrix} $$ since $HT_z = 0$. An easy calculation gives, $ D_PB = F_1^*D_PP $.

These show that each entry of the block operator matrix $V_2 V_3$ is the same as the corresponding entry in the block operator matrix $V_3 V_2$. So,$V_2$ and $V_3$ commute. Hence $(V_1,V_2,V_3)$ is a commuting triple of bounded operators on $\mathcal{K}$.\\
To show that $V_1 = V_2^*V_3 $, we have
\begin{align*}
V_2^*V_3 &= \left[ \begin{array}{c | c c c c c}
B^* & D_PF_1^* & D_PF_1 & 0 & 0 &\cdots \\
\hline
0 & 0 & F_1^* & F_1 & 0 & \cdots \\
0 & 0 & 0 & F_1^* & F_1 & \cdots \\
0 & 0 & 0 & 0 & F_1^* & \cdots \\
0 & 0 & 0 & 0 & 0 & \cdots \\
\cdots & \cdots & \cdots & \cdots & \cdots & \cdots
\end{array}\right]
 \left[\begin{array}{c | c c c c c}
P & 0 & 0 & 0 & 0 & \cdots \\
\hline
0 & 0 & 0 & 0 & 0 & \cdots \\
D_P & 0 & 0 & 0 & 0 & \cdots \\
0 & I & 0 & 0 & 0 & \cdots \\
0 & 0 & I & 0 & 0 & \cdots \\
\cdots & \cdots & \cdots & \cdots & \cdots & \cdots
\end{array}\right] \\
&= \left[ \begin{array}{c | c c c c c}
B^*P + D_PF_1D_P & 0 & 0 & 0 & 0 &\cdots \\
\hline
F_1^*D_P & F_1 & 0 & 0 & 0 & \cdots \\
0 & F_1^* & F_1 & 0 & 0 &\cdots \\
0 & 0 & F_1^* & F_1 & 0 & \cdots \\
0 & 0 & 0 & F_1^* & F_1 & \cdots \\
\cdots & \cdots & \cdots & \cdots & \cdots & \cdots
\end{array} \right] =V_1
\end{align*}
since $B^*P + D_PF_1D_P = A $.
We compute $V_1^*V_1$ to calculate $\lVert V_1 \rVert$. It turns out to be a block diagonal matrix.
\begin{align*}
V_1^*V_1 =  \left[ \begin{array}{c | c c c c}
A^*A+D_PF_1F_1^*D_P & 0 & 0 & 0 &\cdots \\
\hline
0 & F_1^*F_1 + F_1F_1^* & 0 & 0 & \cdots \\
0 & 0 & F_1^*F_1 + F_1F_1^* & 0 &\cdots \\
0 & 0 & 0 & F_1^*F_1 + F_1F_1^* & \cdots \\
0 & 0 & 0 & 0 & \cdots \\
\cdots & \cdots & \cdots & \cdots & \cdots
\end{array} \right]
\end{align*}
Therefore,
$$
\lVert V_1 \rVert ^2 = \lVert V_1^*V_1 \rVert = max \left\lbrace  \lVert A^*A+D_PF_1F_1^*D_P \rVert ,\, \lVert  F_1^*F_1 + F_1F_1^*\rVert \right\rbrace
$$
Using \eqref{F1F2}, we get that
 $$ \lVert A^*A+D_PF_1F_1^*D_P \rVert  =\lVert 0 \oplus \dispmatrix H^*H+HH^* & 0 \\  0 & 0 \\ \rVert =\lVert  F_1^*F_1 + F_1F_1^*\rVert =  \lVert H^*H+HH^* \rVert .$$
The last quantity is $\lVert H_1^*H_1+H_1H_1^* \rVert = \lVert \textmatrix  |\alpha|^2 & 0 \\  0 & |\alpha|^2 \\  \rVert = |\alpha|^2$.
Hence, $\lVert V_1 \rVert = |\alpha| \leq 1$. Since $V_1= V_2^*V_3 $ and $V_3$ is an isometry, we get $V_2=V_1^*V_3$ and $\lVert V_2 \rVert = |\alpha| \leq 1$.\\
\textbf{Conclusion:} The triple $(V_1, V_2, V_3)$ is a commuting triple of bounded operators on $\mathcal{K} $ satisfying
\begin{enumerate}
\item $\mathcal{H}$ is a co-invariant subspace for each $V_j$ (which is clear from their operator matrix representations),

\item $V_3$ is an isometry satisfying $V_1= V_2^* V_3$,

\item $\lVert V_j \rVert = |\alpha| \leq 1$ for $j=1,2$.
\end{enumerate}

Hence by Theorem \ref{thm-e-contraction}, $(V_1, V_2, V_3)$ is a tetrablock isometry. Thus, we have dilated Example \ref{Pal-example} to a tetrablock isometry. It was interesting that we used $T_{z^2}$ on the ortho-complement of $\mathcal{H}$ in $\mathcal{K}$. In section \ref{DilThmSec}, we investigate a class of $(A,B,P)$, each member of which has a tetrablock isometric dilation $(V_1,V_2,V_3)$ with $V_3|_{\mathcal{K}\ominus \mathcal{H}}= T_{z^2}$.

\section{Adjoints} \label{adj}

If $\overline{\mathbb{E}}$ a spectral set for $(A,B,P)$, the Taylor joint spectrum $\sigma (A,B,P)\subseteq \overline{\mathbb{E}}$. The set $\overline{\mathbb{E}}$ is invariant under the complex conjugation i.e., $\overline{\mathbb{E}}^* :=\left\lbrace (\overline{z_1},\overline{z_2},\overline{z_3}):(z_1,z_2,z_3)\in \overline{\mathbb{E}}\right\rbrace = \overline{\mathbb{E}}$. Hence, by the properties of Taylor joint spectrum (see, \cite{Curto})
$$\sigma (A^*, B^* , P^*) = \{ (\overline{z}_1, \overline{z}_2, \overline{z}_3): (z_1, z_2, z_3) \in \sigma (A,B,P)\} \subseteq \overline{\mathbb{E}}.$$

A little work with polynomials then shows that $\overline{\mathbb{E}}$ is a spectral set for $(A^*, B^* , P^*)$ as well. This trick with the polynomials goes through for matrix valued polynomials and hence we have the following proposition. We omit the proof since it is fairly standard.

\begin{proposition}
Let $(A,B,P)$ be a tetrablock contraction such that the tetrablock $\overline{\mathbb{E}}$ is a complete spectral set for $(A,B,P)$. Then $\overline{\mathbb{E}}$ is a complete spectral set for $(A^*,B^*,P^*)$.

\end{proposition}

\begin{remark}
By the proposition above, we see that $\overline{\mathbb{E}}$ is a complete spectral set for $(A^*, B^*, P^*)$ for the triple in Section \ref{sec-1}. By Arveson's theorem, $(A^*, B^*, P^*)$ has a tetrablock isometric dilation. We can actually construct a tetrablock isometric dilation explicitly.  We shall write it down leaving the checking to the reader. Consider $W_1,W_2,W_3$ as follows.
Note that $$ D_{P^*} = \dispmatrix I & 0 \\  0 & I  \\ \oplus \dispmatrix I- T_z T_z^* & 0 \\ 0 & 0\\
 \,\,\, \text{and} \,\,\, \mathcal{D}_{P^*} = \left(\ell ^2(\mathbb{C}^2) \oplus \ell ^2(\mathbb{C}^2)\right) \oplus \left( \text{Ker}(T_z^*)\oplus \{0\} \right).
$$
Define $$ G_1 = 0 \oplus \dispmatrix H^* & 0 \\ 0 & 0 \\
\,\,\text{and}\,\, G_2 =0.
$$
Then
$A^*-BP^* = D_{P^*}G_1 D_{P^*}$ and  $B^* - AP^* = D_{P^*} G_2 D_{P^*}$ by straightforward computations. Define
$$ W_1 :=\left [ \begin{matrix}
A^* & 0 & 0 & 0 & \dots \\
G_1^*D_{P^*} & G_1 & 0 & 0 & \dots \\
0 & G_1^* & G_1 & 0 & \dots \\
0 & 0 & G_1^* & G_1 & \ddots \\
0 & 0 & 0 & G_1^* & \ddots \\
0 & 0 & 0 & 0  & \ddots \\
\vdots & \vdots & \ddots & \ddots & \ddots
\end{matrix}\right ] , \,\,
 W_2 :=
\left [ \begin{matrix}
B^* & 0 & 0 & 0 & \dots \\
G_1D_{P^*} & 0 & 0 & 0 & \dots \\
G_1^* D_{P^*} & G_1 & 0 & 0 & \dots \\
0 & G_1^* & G_1 & 0 & \ddots \\
0 & 0 & G_1^* & G_1 & \ddots \\
0 & 0 & 0 &  G_1^*  & \ddots \\
\vdots & \vdots & \ddots & \ddots & \ddots
\end{matrix}\right ]
$$ and
$$ W_3 := \left [ \begin{matrix}
P^* & 0 & 0 & 0 & \dots\\
0 &	0 & 0 & 0 & \dots \\
D_{P^*} & 0 & 0 & 0 & \dots \\
0 & I & 0 & 0 & \dots \\
0 & 0 & I & 0 & \dots \\
\vdots & \vdots & \ddots & \ddots & \dots
\end{matrix} \right] $$
on $\mathcal{H}\oplus \mathcal{D}_P \oplus \mathcal{D}_P \oplus \cdots$. As before, we can check that $(W_1, W_2,W_3)$ is a tetrablock isometric dilation of $(A^*,B^*,P^*)$.

\end{remark}

\section{A class of tetrablock contractions} \label{DilThmSec}
We investigate the structure of a tetrablock isometric dilation $(V_1,V_2,V_3)$ which has $T_{z^2}$ in the $\mathcal{K}\ominus \mathcal{H}$ part of $V_3$. A tetrablock contraction $(A, B, P)$ gives rise to a pair of operator equations called the {\em fundamental equations} in the unknowns $F_1$ and $F_2$:
$$A - B^*P = D_P F_1 D_P \text{ and } B - A^*P = D_P F_2 D_P.$$
The solutions $F_1$ and $F_2$ exist and are unique. They are called the {\em fundamental operators} of $(A, B, P)$. See  \cite{Bhattacharyya}. For a Hilbert space $\mathcal E$, recall the identification of $\ell ^2(\mathcal E)$ with $H^2_{\mathbb D}(\mathcal E)$. It will be used below.

The beginning of this section warrants a short discussion on Toeplitz operators. The symbol $H^\infty_{\mathbb D} (\mathcal B (\mathcal E))$ will stand for the algebra of operator valued bounded holomorphic functions on the unit disc taking values in  $\mathcal B (\mathcal E)$. Let $\sigma$ be the normalized arc length measure on the unit circle $\mathbb T$. If $\psi$ is in $L^\infty_{\mathbb T}(\mathcal E)$, then $\psi$ induces a linear operator $M_\psi$ on $L^2_{\mathbb T}(\mathcal E)$ by
$$ (M_\psi f)(z) = \psi(z)f(z).$$
Clearly, this is a bounded operator because
\begin{equation} \label{NormOfM} \| M_\psi f \|^2 = \int_{\mathbb T} \|\psi(z)f(z)\|^2 d\sigma(z) \le \|\psi\|_\infty^2 \|f\|^2.\end{equation}

The Hardy space $H^2_{\mathbb D}(\mathcal E)$ of $\mathcal E$-valued holomorphic functions on the unit disc has a natural isometric image inside $L^2_{\mathbb T}(\mathcal E)$ because any Hardy space function of the form $\sum_{n=0}^\infty e_n z^n$ for $e_0, e_1, e_2 , \ldots $ coming from $\mathcal E$ is mapped to the function $\sum_{n=0}^\infty e_n e^{in\theta}$. This is an isometry. Thus, the Hardy space will be treated as a closed subspace of $L^2_{\mathbb T}(\mathcal E)$ and we shall denote the projection onto the Hardy space from $L^2_{\mathbb T}(\mathcal E)$ by $\mathbf{P}$. For $\psi$ as above, define $T_\psi$, the {\em Toeplitz operator with symbol} $\psi$ on the Hardy space by $T_\psi = \mathbf{P} M_\psi |_{H^2_{\mathbb D}(\mathcal E)} $.  We note here that, by virtue of \eqref{NormOfM}, we have
\begin{equation} \label{NormOfT}
\| T_\psi \| \le \|\psi\|_\infty. \end{equation}
For more on Toeplitz operators, we refer the reader to the book \cite{BS}.

Before we state the theorem, it will be good to notice that, in the language of Toeplitz operators, the tetrablock isometry constructed in Section \ref{sec-1} has the following form.
$$V_1= \begin{pmatrix}
  A & 0 \\
  C_1 & T_{\varphi_1}
  \end{pmatrix},\,
   V_2 = \begin{pmatrix}
  B & 0\\
  C_2 & T_{\varphi_2}
  \end{pmatrix}\,\,\text{and}\,\,
  V_3 = \begin{pmatrix}
  P & 0\\
  C_3 & T_{z^2}
  \end{pmatrix}$$
  where
  $$C_1 = (F_1^*D_P, 0, 0, \ldots ), C_2 = (F_1D_P, F_1^*D_P, 0, \ldots ) \text{ and } C_3 = (0, D_p, 0, \ldots , )$$
and
$$\varphi_1(z) =  F_1 + F_1^* z \text{ and } \varphi_2(z) = F_1 z + F_1^* z^2.$$
\begin{thm} \label{main-thm}
Let $(A,B,P)$ be a tetrablock contraction on a Hilbert space $\mathcal{H}$ with fundamental operators $F_1$and $F_2$. On the Hilbert space $\mathcal{K}= \mathcal{H}\oplus \ell ^2(\mathcal{D}_P)$, consider the three bounded operators:
\begin{align}\label{thm-eqn}
  V_1= \begin{pmatrix}
  A & 0 \\
  C_1 & T_{\varphi_1}
  \end{pmatrix},\,
   V_2 = \begin{pmatrix}
  B & 0\\
  C_2 & T_{\varphi_2}
  \end{pmatrix}\,\,\text{and}\,\,
  V_3 = \begin{pmatrix}
  P & 0\\
  C_3 & T_{z^2}
  \end{pmatrix}
\end{align}
for some $C_1,C_2 : \mathcal{H} \to \ell ^2(\mathcal{D}_P)$,
$$
C_3 = \begin{bmatrix}
0 & D_P & 0 & 0 & \dots
\end{bmatrix} ^t : \mathcal{H} \to \ell ^2(\mathcal{D}_P)
$$
and $\varphi_1,\varphi_2 \in L_{\mathbb{T}}^{\infty} \left( \mathcal{B}\left( \mathcal{D}_P \right) \right)$. Then $(V_1,V_2, V_3)$ is a tetrablock isometric dilation of $(A,B,P)$ if and only if $\varphi_1$ and $\varphi_2$ are analytic and are of the form
 \begin{align} \label{phi-1-2}
 \varphi_1(z)= F_1 + \Xi z + F_2^*z^2 ,\,\, \varphi_2(z)= F_2 + \Xi^* z + F_1^* z^2
 \end{align}
 and
 \begin{align}\label{C_1C_2}
 C_1=\begin{bmatrix}
\Xi D_P & F_2^* D_P &0 & 0 &\cdots
\end{bmatrix}^t ,\,\, C_2 = \begin{bmatrix}
\Xi^* D_P & F_1 ^* D_P & 0 & 0 &\cdots
\end{bmatrix}^t
 \end{align}
for some $\Xi \in \mathcal{B}\left( \mathcal{D}_P \right)$ satisfying

\begin{align}
 & \left( \Xi F_1^* - \Xi^*F_2^*\right)D_P P = 0,  \label{item-1} \\
 & \left[ F_2, F_2^* \right] - \left[F_1, F_1^* \right] = \left[\Xi, \Xi^* \right],  \label{item-2} \\
 & \left[F_1, F_2\right] = 0,   \label{item-3} \\
 & \left[ \Xi, F_2\right] =\left[\Xi ^*, F_1\right],  \label{item-4} \\
 & \Xi D_P P = 0 \,\,\text{and}\,\, \Xi^* D_P P =0, \label{item-5}
\end{align}
and
$$
 \lVert \varphi_1 \rVert _{\infty} \bydef \sup_{|z|=1} \lVert \varphi_1(z) \rVert \leq 1 (\text{ equivalently } \| \varphi_2\|_{\infty} \le 1).
$$
\end{thm}

We note that given the form of $V_1,V_2$ and $V_3$ (as in \eqref{thm-eqn}), the subspace $\mathcal{H}$ is co-invariant and
$$
A^* = V_1^*|_{\mathcal{H}} ,\,\, B^*= V_2^*|_{\mathcal{H}},\,\, P^*= V_3^*|_{\mathcal{H}}.
$$
Thus, $(V_1, V_2, V_3)$ is always a dilation of $(A, B, P)$. So, we need to prove that the conditions \eqref{phi-1-2} to \eqref{item-5} are equivalent to $(V_1, V_2, V_3)$ being a tetrablock isometry. This will be achieved by a couple of lemmas.

\begin{lemma}\label{lemma1}
Let $V_1,V_2$ and $V_3$ be as in (\ref{thm-eqn}) with $\varphi_1 , \varphi_2 \in H^\infty_{\mathbb D} (\mathcal B (\mathcal{D}_P)) $. Then $V_1 = V_2^*V_3$ and $V_2 = V_1^*V_3$ if and only if \eqref{phi-1-2} and \eqref{C_1C_2} hold
for some $\Xi \in \mathcal{B}(\mathcal{D}_P)$.
\end{lemma}

 \begin{proof}
 If \eqref{phi-1-2} and \eqref{C_1C_2} are satisfied, then a direct computation of $V_1^*V_3$ and $V_2^*V_3$ proves that $V_1 = V_2^*V_3$ and $V_2 = V_1^*V_3$.

  Conversely, if $V_1,V_2$ and $V_3$ are  as in (\ref{thm-eqn}) and satisfy $V_1= V_2^*V_3$ and $V_2 = V_1^*V_3$, then we have
  $$ \begin{bmatrix}
A & 0 \\
C_1 & T_{\varphi_1}
\end{bmatrix}
= V_1 = V_2^*V_3 = \begin{bmatrix}
B^*P+ C_2^*C_3 & C_2^* T_{z^2}\\
T_{\varphi_2}^*C_3 & T_{\varphi_2}^*T_{z^2}
\end{bmatrix} .$$
Therefore, $$
A-B^*P = C_2^* C_3 ,\,
 T_{\varphi_1} = T_{\varphi_2}^*M_{z^2},\,
C_1 = T_{\varphi_2}^*C_3 , \,\,\text{and}\,\,
 C_2^*T_{z^2}= 0 .$$
 Suppose,
$$\varphi_1(z) = \sum_{j=0}^{\infty} \lambda_j z^j \text{ and } \varphi_2 (z)= \sum_{j=0}^{\infty} \mu_j z^j$$
in $H_{\mathbb{D}}^\infty (\mathcal B\left(\mathcal{D}_P)\right)$, where $\lambda_j, \mu_j $ are bounded operators on $\mathcal{D}_P$. So,
$$  T_{\varphi_1} = T_{\varphi_2}^*T_{z^2} .$$
This implies that
$$ \sum_{j=0}^{\infty} \lambda_j z^j = z^2 \sum_{j=0}^{\infty}\mu_j^* \bar{z}^j  \,\, \text{on the unit circle}\,\, \mathbb{T}. $$
 Comparing the coefficients of $1,z,z^2,\dots$ and $\bar{z}, \bar{z}^2,\dots$ we get,

 $$ \mu_2^* = \lambda_0 ,\,\, \mu_1^* = \lambda_1 , \,\, \mu_0^* = \lambda_2 \,\,\text{and}\,\, \lambda_n =0,\,\, \mu_n = 0 \,\, \forall n \geq 3.$$
Hence, $$\varphi_1 (z) = \lambda_0 + \lambda_1 z + \lambda_2 z^2 \text{ and } \varphi_2(z) = \lambda_2^* + \lambda_1^* z + \lambda_0^* z^2. $$
Also, $C_1 = T_{\varphi_2}^* C_3$ gives,
 $$ C_1 = \begin{bmatrix}
\lambda_2 & \lambda_1 & \lambda_0 & 0 & 0 & \dots\\
0 & \lambda_2 & \lambda_1 & \lambda_0 & 0 & \ddots \\
0 & 0 & \lambda_2 & \lambda_1 & \lambda_0 & \ddots \\
\vdots & \ddots & \ddots & \ddots & \ddots & \ddots
\end{bmatrix}    \begin{bmatrix}
0\\ D_P \\ 0 \\ 0 \\ \vdots
\end{bmatrix}  =  \begin{bmatrix}
\lambda_1 D_P \\
\lambda_2 D_P \\
0\\
0\\
\vdots
\end{bmatrix} . $$
The equation $V_2=V_1^*V_3$ gives $$ \begin{bmatrix}
B & 0 \\
C_2 & T_{\varphi_2}
\end{bmatrix} = \begin{bmatrix}
A^* & C_1^* \\
0 & T_{\varphi_1}^*
\end{bmatrix}  \begin{bmatrix}
P & 0 \\
C_3 & T_{z^2}
\end{bmatrix} = \begin{bmatrix}
 A^*P + C_1^*C_3 & C_1^* T_{z^2} \\
 T_{\varphi_1}^* C_3 & T_{\varphi_1}^* T_{z^2}
\end{bmatrix} $$
which implies the following equations:
$$ B-A^*P = C_1^*C_3 ,\,\, C_2= T_{\varphi_1}^*C_3,\,\, C_1^* T_{z^2} = 0 ,\,\, T_{\varphi_2}= T_{\varphi_1}^*T_{z^2}.$$
Using $C_2= T_{\varphi_1}^*C_3$, we get
$$ C_2 = \begin{bmatrix}
\lambda_0^* & \lambda_1^* & \lambda_2^* & 0 & 0 & \dots\\
0 & \lambda_0^* & \lambda_1^* & \lambda_2^* & 0 & \ddots \\
0 & 0 & \lambda_0^* & \lambda_1^* & \lambda_2^* & \ddots \\
\vdots & \ddots & \ddots & \ddots & \ddots & \ddots
\end{bmatrix}    \begin{bmatrix}
0\\ D_P \\ 0 \\ 0 \\ \vdots
\end{bmatrix} = \begin{bmatrix}
\lambda_1^* D_P \\
\lambda_0 ^* D_P \\
0\\
0\\
\vdots
\end{bmatrix}.$$
Clearly, for $\varphi_1, \varphi_2$, as above $T_{\varphi_2}= T_{\varphi_1}^*T_{z^2} $. Also, it is easy to check that
$$ C_1^* T_{z^2} = 0 \,\, \text{and} \,\, C_2^* T_{z^2} = 0 \,\, \text{for} \,\, C_1, C_2 \,\, \text{as above}. $$
Note that
$$ C_2^*C_3 = D_P \lambda_0 D_P \,\, \text{and}\,\, C_1 ^* C_3 = D_P \lambda_2^* D_P $$

Now we are left with the following two equations :
$$ A-B^*P = C_2^*C_3 = D_P \lambda_0 D_P \text{ and }  B-A^*P = C_1^* C_3 = D_P \lambda_2^* D_P .$$
By uniqueness of the fundamental operators $F_1$ and $F_2$, we get $\lambda_0 = F_1$ and $\lambda_2^* = F_2 $. Set $\lambda_1=\Xi$.
We are done.
 \end{proof}

\begin{lemma}\label{lemma2}
Let $(V_1,V_2,V_3)$ be as in \eqref{thm-eqn}. Suppose \eqref{phi-1-2} and \eqref{C_1C_2} are satisfied. Then $(V_1,V_2,V_3)$ forms a commuting triple if and only if the conditions (\ref{item-1}) to (\ref{item-5}) as in the statement of Theorem \ref{main-thm} hold.
\end{lemma}

\begin{proof}
There are a couple of well-known relationships that the fundamental operators satisfy:
\begin{align}
D_PA = F_1D_P + F_2^*D_P P \text{ and } D_PB = F_2D_P + F_1^*D_PP. \label{item-6}
\end{align}
These were proved in Corollary 4.2 of \cite{Bhattacharyya}. Using these, it is straightforward that the commutativity of $V_1$ and $V_3$ and the commutativity of  $V_2$ and $V_3$ together are equivalent to $(\ref{item-5})$.

Computations show that
$$ V_1V_2 = \left [ \begin{matrix}
AB & 0 & 0 & \cdots \\

\Xi D_PB + F_1\Xi ^*D_P & F_1F_2 & 0 & \cdots \\

F_2^*D_PB + \Xi\Xi ^*D_P + F_1F_1^*D_P  &  \Xi F_2 + F_1\Xi ^* & F_1F_2 & \cdots \\

F_2^*\Xi ^*D_P + \Xi {F_1}^*D_P & F_2^*F_2 + \Xi\Xi^* + F_1F_1^* & \Xi F_2 + F_1\Xi^* & \cdots \\

F_2^*F_1^*D_P & F_2^*\Xi^* + \Xi F_1^* &  F_2^*F_2 + \Xi\Xi^* + F_1F_1^* & \cdots \\

0 & F_2^*F_1^* &  F_2^*\Xi^* + \Xi F_1^* & \cdots \\

\cdots & \cdots & \cdots & \cdots &
\end{matrix} \right ]
$$
and
$$
V_2V_1 = \left[\begin{matrix}

BA & 0 & 0 &  \cdots \\

\Xi ^*D_PA + F_2\Xi D_P & F_2F_1 & 0 & \cdots \\

F_1^*D_P A + \Xi ^*\Xi D_P + F_2F_2^*D_P & \Xi ^*F_1+F_2\Xi  & F_2F_1 & \cdots \\

F_1^*\Xi D_P + \Xi ^*F_2^*D_P & F_1^*F_1 + \Xi ^*\Xi + F_2F_2^* & \Xi ^*F_1+F_2\Xi & \cdots \\

F_1^*F_2^*D_P & F_1^*\Xi+ \Xi ^*F_2^* & F_1^*F_1 + \Xi ^*\Xi +F_2F_2^* &\cdots \\

0 & F_1^*F_2^* & F_1^*\Xi + \Xi ^*F_2^* & \cdots \\

\cdots & \cdots & \cdots & \cdots

\end{matrix}\right ]
.$$
 The commutativity of $V_1 $ and $V_2$ is equivalent to the equations from (\ref{item-1}) to (\ref{item-4}) by virtue of the relations \eqref{item-6}.
\end{proof}
\begin{proof}{(of Theorem \ref{main-thm})}

If the conditions \eqref{phi-1-2} to \eqref{item-5} hold, then we first get by Lemma \ref{lemma1} that $V_1= V_2^*V_3$ and $V_2 = V_1^*V_3$. Also by Lemma \ref{lemma2}, $(V_1,V_2,V_3)$ forms a commuting triple. Clearly $V_3$ is an isometry. Since $\varphi_1$ and $\varphi_2$ satisfy $\varphi_1(z) = \varphi_2(z)^* z^2$ for every $z$ on the unit circle, $\| \varphi_1(z) \| = \| \varphi_2(z) \|$ for every $z$ on the unit circle and hence $\lVert \varphi_1 \rVert_{\infty} = \lVert \varphi_2 \rVert_{\infty}$.
So, for $j=1,2$, $ \lVert \varphi_j \rVert_{\infty} \leq 1$ which give $r\left( T_{\varphi_j}\right) \leq 1 $ where $r$ stands for spectral radius. Since $V_j$'s are lower triangular as in (\ref{thm-eqn}),
$$
\sigma (V_1) \subseteq \sigma (A) \cup \sigma (T_{\varphi_1}) \,\, and \,\, \sigma (V_2) \subseteq \sigma (B) \cup \sigma (T_{\varphi_2}).
$$
These containments of spectrum along with the fact that $A,B$ are contractions imply $r(V_j)\leq 1$ for $j=1,2$. Hence by Theorem \ref{thm-e-contraction}, $(V_1,V_2,V_3)$ is a tetrablock isometry.

For the other direction, we first recall a basic fact about Toeplitz operators.
Every $\varphi$ in $L^\infty_{\mathbb T}(\mathcal B (\mathcal E))$ has the form $\sum_{n=-\infty}^\infty \Phi_n e^{in\theta}$ for some $\Phi_n $ in $\mathcal B(\mathcal E)$. The result says that the block operator matrix of $T_\varphi$ with respect to the standard basis of $H^2_{\mathbb D}(\mathcal E)$ is of the form
\begin{align} \label{TphiMatrix}
T_{\varphi}= \begin{bmatrix}
\Phi_0 & \Phi_{-1} & \Phi_{-2} & \Phi_{-3} & \cdots \\
\Phi_1 & \Phi_0 & \Phi_{-1} & \Phi_{-2}&\cdots \\
\Phi_2 & \Phi_1 & \Phi_0 & \Phi_{-1} & \cdots \\
\Phi_3 & \Phi_{2} & \Phi_{1} & \Phi_0 &\cdots\\
\cdots & \cdots & \cdots & \cdots & \cdots
\end{bmatrix}.\end{align}
Note that the $\Phi_n$ are $0$ for all $n < 0$ if and only if $\varphi \in H^\infty_{\mathbb D} \left(\mathcal{B}(\mathcal E)\right) $.

If $(V_1,V_2,V_3)$ is a  tetrablock isometry, then it is a commuting triple.
Commutativity of $V_1$ and $V_3$ is the same as
$$
\begin{bmatrix}
AP & 0\\
C_1P + T_{\varphi_1} C_3 & T_{\varphi_1}T_{z^2}
\end{bmatrix}
= \begin{bmatrix}
PA & 0 \\
C_3A+T_{z^2}C_1 & T_{z^2}T_{\varphi_1}
\end{bmatrix}
$$
which implies that $T_{\varphi_1}$ commutes with $T_{z^2}$. Similarly, commutativity of $V_2$ and $V_3$ implies that $T_{\varphi_2}$ commutes with $T_{z^2}$. So, we need the general structure of an element of the commutator algebra of $T_{z^2}$. By a straightforward computation, we get that every such element has the form
\begin{align}\label{commutant}
T= \begin{bmatrix}
A_{11} & A_{12} & 0 & 0 & 0 & \cdots \\
A_{21} & A_{22} & 0 & 0 & 0 & \cdots \\
A_{31} & A_{32} & A_{11} & A_{12} & 0 & \cdots \\
A_{41} & A_{42} & A_{21} & A_{22} & 0 & \cdots \\
A_{51} & A_{52} & A_{31} & A_{32} & A_{11} & \cdots\\
\cdots & \cdots & \cdots & \cdots & \cdots & \cdots
\end{bmatrix}
\end{align}
in the standard basis of $H^2_{\mathbb D}(\mathcal D_P)$. Looking at \eqref{TphiMatrix} and \eqref{commutant}, we realize that a Toeplitz operator $T_\varphi$ in the commutant algebra of $T_{z^2}$ necessarily has an analytic symbol $\varphi$. Thus, $\varphi_1$ and $\varphi_2$ are in $H^\infty_{\mathbb D}(\mathcal B (\mathcal D_P))$.

We now invoke Theorem \ref{thm-e-contraction} which in particular tells us that $V_1 = V_2^* V_3$ and $V_2 = V_1^* V_3$. So, we can apply Lemma \ref{lemma1} which gives us  \eqref{phi-1-2} and \eqref{C_1C_2}. As soon as we get \eqref{phi-1-2} and \eqref{C_1C_2}, we apply Lemma \ref{lemma2}. We get (\ref{item-1}) to (\ref{item-5}). The norm condition gives $\lVert \varphi_1 \rVert_{\infty} \leq 1$.
\end{proof}

We comment that $\Xi = 0$ is a choice and for this choice the equations \eqref{item-1} - \eqref{item-5} reduce to just two conditions: $\left[F_1, F_2\right] = 0$ and $\left[F_1^*, F_1 \right]= \left[ F_2^*, F_2\right]$. Recall that the minimal dilation under these two conditions was found in \cite{Bhattacharyya}. The one found here is not minimal under these conditions. What worked here is the Toeplitz structure which quickly led to the symbols to be in $ H_{\mathbb{D}}^{\infty} \left( \mathcal{B}\left( \mathcal{D}_P \right) \right)$ which is a proper subalgebra of the commutant algebra of $T_{z^2}$ acting on $H^2_{\mathbb D}\left( \mathcal{D}_P \right)$.

\textbf{Acknowledgement :} The first author's work is supported by J C Bose National Fellowship JCB/2021/000041. The second author's work is supported partly by the National Board for Higher Mathematics Fellowship NBHM-0203/11/2020/R\& D-II/10167 and partly by the Prime Minister's Research Fellowship PMRF-21-1274.03.




\begin{thebibliography}{99}

\bibitem{agler} J. Agler, Rational dilation on an annulus, \textit{Annals of Math.} \textbf{121} (1985), 537-563.

\bibitem{ahr} J. Agler, J. Harland and B. J. Raphael, Classical function theory,
operator dilation theory, and machine computation on multiply-connected
domains, \textit{Mem. Amer. Math. Soc.}  \textbf{191} (2008), no. 892.

\bibitem{ay-jfa} J. Agler and N. J. Young, A commutant lifting theorem for a domain in $\mathbb{C}^2$ and spectral
interpolation, \textit{J. Funct. Anal.} \textbf{161} (1999), 452-477.

\bibitem{Ando}
 T. Ando, On a pair of commutative contractions,\textit{Acta Sci. Math.} \textbf{24} (1963), 88-90.

\bibitem{sub2} W. Arveson, Subalgebras of $C^*$-algebras II,
\textit{Acta Math.}, \textbf{128} (1972), no. 3-4, 271–308.

\bibitem{Ball-Sau}
 J. A. Ball, H. Sau, Rational dilation of tetrablock contractions revisited,\textit{J. Funct. Anal.} \textbf{278} (2020), no. 1, 108275, 14 pp.

 \bibitem{Bhattacharyya}
T. Bhattacharyya, The tetrablock as a spectral set, \textit{Indiana Univ. Math. J.} \textbf{63} (2014), no. 6, 1601–1629.

\bibitem{BS} A. Böttcher and  B. Silbermann, Analysis of Toeplitz Operators, 2nd edition,
Springer-Verlag, Berlin, Heidelberg, New York 2006.

\bibitem{Curto} R. E. Curto, Applications of several complex variables to multiparameter spectral theory,Surveys of some recent results in operator theory, Vol. II, 25–90,\textit{Pitman Res. Notes Math. Ser.}, \textbf{192}, Longman Sci. Tech., Harlow, 1988.

\bibitem{DM} M. A. Dritschel and S. McCullough, The failure of rational dilation on a triply connected domain,
\textit{ J. Amer. Math. Soc.}  \textbf{18} (2005), 873-918.


\bibitem{Du-Jin} H. K. Du and P. Jin, Perturbation of spectrums of $2\times 2$ operator matrices, \textit{Proc. Amer.Math. Soc.} \textbf{121}(1994), no. 3, 761-766.

\bibitem{vN}
von Neumann, Eine Spektraltheorie für allgemeine Operatoren eines unitären Raumes (German),\textit{Math. Nachr.} \textbf{4} (1951), 258–281.

\bibitem{Pal}
S. Pal, The failure of rational dilation on the tetrablock, \textit{J. Funct. Anal.} \textbf{269} (2015), no. 7, 1903–1924.

\bibitem{Paulse-JFA} V. I. Paulsen, Every completely polynomially bounded operator is similar to a contraction, \textit{J. Funct. Anal.} \textbf{55} (1984), 1–17.

\bibitem{Pisier-os} G. Pisier, Introduction to Operator Space Theory, \textit{Cambridge University Press}, 2003.


\bibitem{Var}N. Th. Varopoulos, On an inequality of von Neumann and an application of the metric theory of tensor products to operators theory, \textit{J. Funct. Anal.} \textbf{16} (1974), 83–100.

























\end{thebibliography}


\end{document}